\documentclass{amsart}
\usepackage{amsfonts,amssymb,amsmath,amsthm}
\usepackage[latin1]{inputenc}
\usepackage{multicol}
\usepackage[all]{xy}
\usepackage{multirow}

\newcommand{\ballontop}[1]{\text{\accent'27 #1}}

\newcommand{\gr}[1]{{{\mbox{\scriptsize $(#1)$}}}}

\newtheorem{thm}{Theorem}[section]
\newtheorem{theorem}[thm]{Theorem}

\newtheorem{proposition}[thm]{Proposition}

\newtheorem{lemma}[thm]{Lemma}

\theoremstyle{definition}
\newtheorem{defn}[thm]{Definition}
\newtheorem{definition}[thm]{Definition}

\newtheorem{rem}[thm]{Remark}

\newtheorem{example}[thm]{Example}

\numberwithin{equation}{section}

\title{Desingularization of Legendrian Surfaces - The Quasi-Ordinary case}

\author{Ant\'onio Ara\'ujo, Jo\~ao Cabral and Orlando Neto}

\begin{document}
\maketitle



\begin{abstract} In this paper we prove a desingularization theorem for Legendrian surfaces that are the conormal of a quasi-ordinary hypersurface. 
\end{abstract}

\section{Introduction}

\noindent
The main purpose of this paper is to prove a desingularization theorem for Lagrangean surfaces of a contact manifold of dimension 5. For the moment we limit ourselves to considering the quasi-ordinary case: the case when the Lagrangean variety is the conormal of a quasi-ordinary surface. Our proof depends on recent work on the computation of the limits of tangents of a quasi-ordinary surface \cite{LIMITS}. At the moment we do not have a systematic way to compute the limits of tangents of a general class of hypersurfaces. Once this problem is solved, it should not be very hard to generalize our main result to an arbitrary Lagrangean variety.

\noindent
Our paper generalizes previous work on the desingularization of Lagrangean curves (see \cite{neto3}). Lipman proved a desingularization theorem for quasi-ordinary surfaces (see \cite{Lipmanphd}). Ban and MacEwan \cite{mcewan} showed that Lipman's algorithm produces an embedded desingularization, which coincides with the algorithm of Bierston and Milman \cite{milman}.  

\noindent
We show that when a quasi-ordinary surface $S$ has trivial limits of tangents and $L$ is an admissible center for $S$, the conormal of the blow up of $S$ along $L$ equals the blow up of the conormal of $S$ along the conormal of $L$. 

\noindent 
We recall that each Lagrangean variety of $\mathbb{P}^*X$ is the conormal of its projection on $X$. Moreover, each Lagrangean variety is isomorphic to the conormal of a hypersurface with trivial limits of tangents (see \cite{kk}). Hence we can apply the procedure to each germ of Lagrangean surface.

\noindent
One of the main motivations of this work is its application to the desingularization of certain classes  of holonomic systems of partial differential equations (see \cite{neto5}).


\section{Logarithmic contact manifolds}

\noindent
All manifolds considered in this paper are complex analytic manifolds. 

\noindent
A subset $Y$ of a manifold $X$ is called a divisor with normal crossings at $o\in X$ if there is an open neighborhood $U$ of $o$, a system of local coordinates $(x_1,\ldots,x_n)$ and a nonnegative integer $\nu$ such that $x_i(o)=0, i=1,\ldots ,n$, and
\begin{equation}\label{normal_crossings_coord}
Y \cap U=\{x_1 \cdots x_{\nu}=0\}.
\end{equation} 
\noindent We call $\nu$ the index of $Y$ at $o$. We say that $Y$ is a \emph{divisor with normal crossings} if $Y$ is a divisor with normal crossings at each point of $X$. 

\noindent
Let $Y_1,\ldots,Y_\nu$ be the irreducible components of  $Y$ in a neighboorhood of  $o$.
Let $1\le \mu\le \nu$. Set $Z=Y_1\cap\cdots\cap Y_\mu$. We call $(Y_{\mu+1}\cup \cdots \cup Y_\nu)\cap Z$ the \em normal crossings divisor induced in $Z$ by $Y$. \em 

\noindent  
Let $j:X\setminus Y \hookrightarrow X$ be the open inclusion. 
Let $\mathcal{O}_X$ denote the sheaf of holomorphic functions of $X$. 
Let $\Omega_X^1$ denote the sheaf of differential forms of degree $1$ on $X$.
Let $\Omega_X^{1}\langle Y \rangle$ be the smallest subsheaf of 
$j_{*}\Omega^{1}_{X \setminus Y}$ that contains $\Omega_X^{1}$ and $d\log f$ for
each holomorphic function $f$ such that $f^{-1}(0)\subset Y$.
Set $\Omega_X^{p}\langle Y \rangle=\wedge^p\Omega_X^{1}\langle Y\rangle$.
 The local sections of $\Omega_X^{*}\langle Y \rangle$ are called logarithmic differential forms with poles along $Y$.

\noindent
Let $Z$ be a smooth irreducible component of $Y$. We can associate to $\alpha \in \Omega_X^1\langle Y\rangle$ an holomorphic function 
$Res_Z \alpha \in \mathcal{O}_Z$, the \emph{Poincar\'{e} residue of $\alpha$ along $Z$}.  
Assume that we are in the situation of (\ref{normal_crossings_coord}),
$\alpha|_U=\sum_{i=1}^{\nu}\alpha_i {d x_i}/{x_i}+\sum_{i=\nu+1}^n\alpha_idx_i$
and $Z \cap U=\{x_j=0\}$, where $1 \leq j \leq \nu$. Then 
$Res_Z\alpha|_{U\cap Z}=\alpha_j|_{U\cap Z}.$
Let $W$ be the intersection of the smooth irreducible components $Y_1, \ldots, Y_{\mu}$ of $Y$. 
We call \emph{residual submanifold of $X$ along $W$} to the set of points $o \in W$ such that  Res$_{Y_i}(\theta)$  vanishes at $o$ for $1 \leq i \leq \mu$. 
We will denote the residual submanifold of $X$ along $W$ by $R_WX$

\vskip .1in
\noindent 
A group action $\alpha: \mathbb{C}^{*}\times X \rightarrow X$ on manifold $X$ is called a free group action of $\mathbb{C}^{*}$ on $X$ if, for each $x \in X$, the isotropy subgroup $\{t \in \mathbb{C}^{*}:\alpha(t,x)=x\}$ equals $\{1\}$. A manifold $X$ with a free froup action $\alpha$ of $\mathbb{C}^{*}$ is called a homogeneous manifold. We associate to each free group action $\alpha$ of $\mathbb{C}^{*}$ on $X$ a vector field $\rho$, the \emph{Euler vector field of $\alpha$}, given by
 $$\rho f=\frac{\partial}{\partial t}\alpha_x^{*}f|_{t=1}, f \in \mathcal{O}_X,$$
where $\alpha_x(t)=\alpha (t,x)$. 
 Given homogeneous  manifolds $(X_1,\alpha_1)$ and $(X_2,\alpha_2)$, 
 a holomorphic map $\varphi:X_1 \rightarrow X_2$ is called \emph{homogeneous} if 
 $\alpha_{2,t}\varphi=\varphi\alpha_{1,t},$
 for any $t \in \mathbb{C}^{*}.$

\noindent 
Let us recall some definitions and some results introduced in \cite{neto4}.

\begin{definition}
Let $X$ be a  manifold of dimension $2n$ endowed with a free group action $\alpha$. 
Let $Y$  be a divisor with normal crossings of $X$. 
If $\sigma$ is a locally exact section of $\Omega_X^2\langle Y \rangle$, 
$\alpha^*_t\sigma=t\sigma$ for each $t\in\mathbb C^*$ and 
$\sigma^n$ is a generator of $\Omega_X^{2n}\langle Y \rangle$
we say that $\sigma$
is a \emph{homogeneous logarithmic symplectic form with poles along $Y$} and $(X,\sigma)$ a 
\emph{homogeneous logarithmic symplectic manifold with poles along $Y$}. 

\noindent  
Let $(X_1,\sigma_1), (X_2,\sigma_2)$ be homogeneous symplectic manifolds. 
A homogeneous map $\varphi:X_1 \rightarrow X_2$ is a \em homogeneous symplectic transformation \em 
if $\varphi^*\sigma_2=\sigma_1$.
\end{definition}

\noindent
If $Y$ is the empty set we get the usual definition of homogeneous symplectic manifold.

\noindent 
Given a homogeneous logarithmic symplectic manifold $(X,\sigma)$ we call 
$\theta= \iota (\rho)\sigma$
the \emph{canonical 1-form} of $(X,\sigma)$, where $\iota(\rho)\sigma$ is the contraction of $\rho$ and $\sigma$. 
Notice that $\sigma=d\theta$.

\noindent
Given a vector bundle $E$ over a manifold $M$ we denote by $\ballontop{E}$ the  manifold $E\setminus M$, where we identify $M$ with the image of the zero section of $E$.

\begin{example}
Let $M$ be a  manifold and $N$ a divisor with normal crossings of $M$. Let
$\pi :T^{*}\langle M / N \rangle \rightarrow M$
be the vector bundle with sheaf of sections $\Omega_M^1\langle N \rangle$. 
We will call $T^{*}\langle M / N\rangle$ the \emph{logarithmic cotangent bundle of $M$ along $N$}. 

\noindent
The manifold $\ballontop{T}^{*}\langle M/N  \rangle $ has a canonical structure of logarithmic symplectic manifold with poles along $\pi^{-1}(N)$. 
Actually, there is a canonical section $\theta$ of $\Omega^1_{T^{*}\langle M / N \rangle}\langle \pi^{-1}(N)\rangle$, the canonical \emph{1-form} of $T^{*}\langle M / N \rangle$.
Given an integer $\nu$ and a system of local coordinates $(x_1, \ldots, x_n)$ on an open set $U$ of $X$ verifying ($\ref{normal_crossings_coord}$) there is one and only one family of holomorphic functions $\xi_i, 1 \leq i \leq n$, defined on $\pi^{-1}(U)$ such that $\theta|_{\pi^{-1}(U)}$ equals
\begin{equation}\label{teta}
\sum_{i=1}^{\nu}\xi_i\frac{dx_i}{x_i}+\sum_{i=\nu +1}^{n}\xi_i dx_i.
\end{equation}
Moreover, the 2-form $\sigma=d\theta$  is a homogeneous symplectic form with poles along $\pi^{-1}(Y)$.
\end{example}

\noindent 
A homogeneous logarithmic symplectic manifold is locally isomorphic to $\mathring{T}^{*}\langle M/N \rangle$ in the category of homogeneous symplectic manifolds:

\begin{theorem}
Let $\sigma$ be a homogeneous logarithmic symplectic form on a  manifold $X$ with poles along a divisor with normal crossings $Y$. 
Given $o \in X$ let $\nu$ be the number of irreducible components of $Y$ at $o$. Then there is a system of local coordinates 
\begin{equation}\label{SCOORDINATES}
(x_1,\ldots , x_n,\xi_1, \ldots, \xi_n )
\end{equation}

\noindent
on an open conic neighbourhood $U$ of $o$ such that $Y\cap U=\{x_1\cdots x_{\nu}=0\}$, $x_1, \ldots , x_n$ are homogeneous of degree $0$, 
$\xi_1, \ldots, \xi_n$ are homogeneous of degree $1$ and $\sigma |_U$ equals
\begin{equation}\label{SNORMALFORM}
\sum_{i=1}^{\nu}d\xi_i\frac{dx_i}{x_i}+\sum_{i=\nu +1}^{n}d\xi_i dx_i.
\end{equation}
\end{theorem}

\noindent
A symplectic form $\sigma$ on a manifold $X$ defines on $X$ a Poisson bracket $\{\cdot,\cdot\}$.
Moreover, we can recover $\sigma$ from the Poisson bracket. A submanifold $W$ of $X$ is called \em involutive \em [\em invariant\em] if $\{I_W, I_W\}\subset I_W$ [$\{I_W, \mathcal{O}_X\}\subset I_W]$.

\noindent
Given coordinates (\ref{SCOORDINATES}) on a conic neighbourhood $U$ of $o$ such that $\sigma|_U$ equals (\ref{SNORMALFORM}), 
we have that $\iota (\rho)\sigma |_U$ equals (\ref{teta}) and $\{f,g\}$ equals

\begin{equation}\label{POISSON}
\sum_{i=1}^{\nu}x_i\left(\frac{\partial f}{\partial \xi_i}\frac{\partial g}{\partial x_i}-\frac{\partial f}{\partial x_i}\frac{\partial g}{\partial \xi_i}\right)+\sum_{i=\nu +1}^{n}\left(\frac{\partial f}{\partial \xi_i}\frac{\partial g}{\partial x_i}-\frac{\partial f}{\partial x_i}\frac{\partial g}{\partial \xi_i}\right)
\end{equation}
for each $f,g\in\mathcal O_X(U)$.

\begin{proposition}
Let $X$ be an homogeneous logarithmic symplectic manifold with poles along a smooth divisor $Y$. 
Let $W$ be the intersection of the smooth irreducible components $Y_1, \ldots, Y_{\mu}$ of $Y$.  
Then:

\noindent
$(i)$
{\em $X$, $R_WX$ are involutive submanifolds of $X$.\em}

\noindent
$(ii)$
{\em The manifold $R_WX$ has a canonical structure of homogeneous symplectic manifold with poles along the divisor induced in $W$ by $Y$.\em}

 \end{proposition}
 \begin{proof}
 Let $o \in W$. There is a system of symplectic coordinates (\ref{SCOORDINATES}) 
 on a conic open set $U$ that contains $o$ such that 
 $\theta|_U$ equals (\ref{teta}) and $W\cap U=\{x_1=\cdots=x_{\mu}=0\}$. 
 Hence $R_WX\cap U=\{x_1=\cdots=x_{\mu}=\xi_1=\cdots=\xi_{\mu}=0\}.$
 The restriction to $R_WX\cap U$ of the Poisson bracket of $X$ is given by
 $$
 \{f,g\}=\sum_{i=\mu+1}^{\nu}x_i\left(\frac{\partial f}{\partial \xi_i}\frac{\partial g}{\partial x_i}
 -\frac{\partial f}{\partial x_i}\frac{\partial g}{\partial \xi_i}\right)+
 \sum_{i=\nu +1}^{n}\left(\frac{\partial f}{\partial \xi_i}\frac{\partial g}{\partial x_i}-
 \frac{\partial f}{\partial x_i}\frac{\partial g}{\partial \xi_i}\right).
 $$
Hence $R_WX\cap U$ is endowed with a 1-form  
 $\sum_{i=\mu+1}^{\nu}\xi_i{dx_i}/{x_i}+\sum_{i=\nu+1}^n\xi_idx_i$.
 \end{proof}

 \begin{definition}
 Let $X$ be a  manifold of dimension $2n+1$ and $Y$ a divisor with normal crossings of $X$. A local section $\omega$ of $\Omega^1_X\langle Y \rangle$ is called a \emph{logarithmic contact form with poles along $Y$} if $\omega (d\omega)^n$ is a local generator of $\Omega_X^{2n+1}\langle Y \rangle$.
 
 \noindent We say that a locally free sub $\mathcal{O}_X$-module $\mathcal{L}$ of $\Omega^1_X\langle Y \rangle$ is a \emph{logarithmic contact structure on $X$ with poles along $Y$} if it is locally generated by a logarithmic contact forms with poles along $Y$. We say that a  manifold with a logarithmic contact structure with poles along a divisor with normal crossings $Y$ is a \emph{logarithmic contact manifold with poles along $Y$}. We call $Y$ the  \emph{set of poles of the logarithmic contact manifold $(X,\mathcal{L})$.}
 
 \noindent Let $(X_1,\mathcal{L}_1)$, $(X_2,\mathcal{L}_2)$ be logarithmic contact manifolds. We say that a holomorphic map $\varphi:X_1 \rightarrow X_2$ is a \emph{contact transformation} if $\varphi^*\omega$ is a local generator of $\mathcal{L}_1$ for each local generator $\omega$ of $\mathcal{L}_2$  .
 \end{definition}

\noindent Let $X$ be a homogeneous logarithmic symplectic manifold. 
Let $\theta$ be the canonical 1-form of $X$ and let $Y$ be the set of poles of $X$. Let $X_{*}$ be the quotient of $X$ by its $\mathbb{C}^{*}$ action. Then $X_{*}$ is a  manifold and the canonical epimorphism $\gamma:X \rightarrow X_{*}$ is a $\mathbb{C^{*}}$-bundle. Put $Y_{*}=\gamma (Y)$. 
Let $\mathcal{L}$ be the sub $\mathcal{O}_{X_{*}}$-module of $\Omega_{X_{*}}^1\langle Y \rangle$ generated by the logarithmic differential forms $s^{*}\theta$, where $s$ is a holomorphic section of $\gamma$. 
Then $\mathcal{L}_{*}$ is a structure of logarithmic contact manifold with poles along $Y_{*}$.

\noindent
We constructed on this way a functor \em contactification \em $X\mapsto X_*$ 
from the category $\mathcal S$ of homogeneous logarithmic symplectic manifolds into the category $\mathcal C$ of logarithmic contact manifolds. 
Let $\mathcal S'$ [$\mathcal C'$] be the subcategory of $\mathcal S$ [$\mathcal C$]
such that its morphisms are locally injective.
The functor contactification defines an equivalence of categories from $\mathcal S'$ onto $\mathcal C'$.

\noindent
Theorem \ref{inversecommutes} constructs an example of a 
morphism of homogeneous logarithmic symplectic manifolds that is not locally injective.
Theorem \ref{diagteo} constructs an example of a 
morphism of logarithmic contact manifolds that is not locally injective.
Moreover, this last morphism does not come from a 
morphism of homogeneous logarithmic symplectic manifolds through the functor $X\mapsto X_*$.
This is the main reason why there is no equivalence of categories between 
$\mathcal S$ and $\mathcal C$ and we have to use logaritmic contact manifolds in this paper.
This phenomena has no equivalent in the classic, non logarithmic case.

\noindent 
It is common to use the coordinates of $\mathbb{C}^n$ when dealing with the projective space $\mathbb C\mathbb{P}^{n-1}$. 
 We will use symplectic coordinates when dealing with logarithmic contact manifolds within the same spirit. 
 In particular we do not feel the need to define concepts like involutivity or residual set in the contact context.

\noindent 
The \emph{ projective logarithmic cotangent bundle of $M$ with poles along $N$} \allowbreak ${\mathbb{ P}}^*\langle M/N \rangle=(\mathring{T}^*\langle M/N \rangle)_*$ has a canonical structure of logarithmic contact manifold.


\section{Legendrian Varieties}
Let $(X,\mathcal L)$ be a contact manifold of dimension $2n+1$. An
analytic subset $\Gamma$ of $X$ is a Legendrian variety of $X$ if
it verifies  the following three conditions: $\Gamma$ has
dimension $n$, $\Gamma$ is involutive and the restriction to the
regular part of $\Gamma$ of a
  local generator of $\mathcal L$ vanishes.

\noindent Each two of these three conditions imply the remaining one.

\noindent Given a  manifold $M$ and an irreducible analytic subset $S$ of $M$ there is 
one and only one Legendrian variety $\mathbb{P}^*_SM$ of $\mathbb{P}^*M$ such that 
$\pi (\mathbb{P}^*_SM)=S$.
The analytic set  $\mathbb{P}^*_SM$ is called the \em conormal \em of $S$ (see for
instance \cite{ka3}). If $S$ has irreducible components $S_i, i \in I$, the conormal
$\mathbb{P}^*_{S}M$ of $S$ equals $\cup_{i \in I}\mathbb{P}^*_{S_j}M.$

\noindent Let us introduce  stratified versions of the definitions above.

\begin{defn}\label{defleg}
Let $X$ be a logaritmic contact manifold of dimension $2n+1$ with
set of poles $Y$. An analytic subset $\Gamma$ of $X$  is called a
Legendrian variety of $X$ if the following conditions are verified:

\noindent
(i)  $\Gamma$ is involutive and  $\Gamma\setminus Y$ is a Legendrian variety of $X\setminus Y$.

\noindent
(ii) Let $o$ be a point of $\Gamma$. 
Let $\Gamma'$ [$Y'$] be an irreducible component of the germ of $\Gamma$ [$Y$] at $o$. If $\Gamma'\subset Y'$,  $\Gamma'\subset R_{Y'}X$.

\noindent
(iii) If $\Gamma''$ is an irreducible component of $\Gamma\cap Y'$ and $\Gamma''\not\subset Y_{sing}$, $\Gamma''$ is a Legendrian subvariety of $R_{Y'}X$.
\end{defn}

\noindent
Let $N$ be a divisor with normal crossings
of a manifold $M$.

\begin{rem}\label{MNQRbis}
 Let $\Gamma$ be a Legendrian variety of $\mathbb{P}^*\langle
M/N\rangle$. Let $Q$ be an irreducible component of $N$. Let $R$ be the divisor
with normal crossings induced in $Q$ by $N$. If $\Gamma$ is
contained in $\pi^{-1}(Q)$, it follows from condition (ii) of
Definition \ref{defleg} that $\Gamma$ is contained in $\mathbb{P}^*\langle
Q/R\rangle$.
\end{rem}

\noindent An analytic subset $S$ of $M$  is a \em natural analytic subset \em of $(M,N)$ if no germ of
$S$ is an intersection of irreducible components of a germ of $N$. 

\noindent 
A Legendrian variety of a logarithmic contact manifold $X$ with
poles along $Y$ is a natural analytic subset of $(X,Y)$.

\begin{defn}
Let $S$ be a natural irreducible subset of $(M,N)$. Let $Q$ be the
intersection of the irreducible components of $N$ that contain $S$.
Let $R$ be the divisor with normal crossings induced in $Q$ by $N$.
We call \em conormal \em
of $S$ (relative to $N$) to the closure $\mathbb{P}^*_S\langle Q/R\rangle$  of the conormal
of the analytic subset $S\setminus R$ of $Q\setminus R$ in
$\mathbb{P}^*\langle Q/R\rangle$.

\noindent Let $S$ be a natural analytic subset of $(M,N)$. We call \em
conormal \em of $S$ to the union $\mathbb{P}^*_S\langle M/N\rangle$ of the
conormals of its irreducible components.
\end{defn}

\begin{example}
Set $M=\mathbb C^4$, $N=\{ x_1x_2x_3=0\}$, $S=\{x_1=x_2=x_4=0\}$.
Hence $Q=\{ x_1=x_2=0\}$, $R=\{x_1=x_2=x_3=0\}$ and 
$\mathbb P^*\langle Q/R\rangle=\{ x_1=x_2=\xi_1=\xi_2=0 \}$
is endowed with the canonical $1$-form $\xi_3dx_3/x_3+\xi_4dx_4$. 
Therefore $\mathbb P^*_S\langle M/N\rangle =\{  x_1=x_2=\xi_1=\xi_2=x_4=\xi_3=0\}$.
\end{example}

\begin{thm}\label{conleg}
The conormal of a natural analytic set is a Legendrian variety.
\end{thm}

\begin{proof}
Let $S$ be a germ of a natural analytic subset of $(M,N)$. Set $\Gamma=\mathbb{P}^*_S\langle M / N \rangle$. We can
assume that $S$ is irreducible and that $M$ is the intersection of the irreducible components of $N$
that contain $S$. 
The intersection
of $\Gamma$ with $\pi^{-1}(M\setminus N)$ is  the
Legendrian variety $\mathbb{P}^*_{S\setminus N}( M\setminus N)$ of the
contact manifold $\mathbb{P}^*( M\setminus N)$. 
Since $\Gamma$ is the closure of $\mathbb{P}^*_{S\setminus N}(M\setminus N)$, $\Gamma$ is involutive. 
Hence condition (i) is verified. 
Condition (ii) follows from the definition of conormal variety. 

\noindent Let us prove condition (iii) by induction in the dimension of $M$.
Condition (iii) is trivial if dim$M=1$. 
Let $Z$ be an irreducible component of $\pi^{-1}(N)$.
The set $Q=\pi(Z)$ is an irreducible component of $N$.
Since $\pi^{-1}(N)$ is invariant, $Z$ is invariant. 
Let $R$ be the divisor induced in $Q$ by $N$. 
Let $\Gamma_0$ be an irreducible component of $\Gamma \cap Z$ that is not contained in the singular locus of $\pi^{-1}(N)$.
Let us show that
\begin{equation}\label{contido}
  \Gamma_0\subset  \mathbb{P}^*\langle Q/R\rangle.
\end{equation}
It is enough to show that $\gamma^{-1}(\Gamma_0)$ is contained in  $T^*\langle Q/R\rangle$. 
Let $o \in \gamma^{-1}(\Gamma_0\setminus \pi^{-1}(N)_{sing})$.  
 There is an open conic
neighborhood $U$ of $o$ and a system of local coordinates
(\ref{SCOORDINATES}) on $U$ such that
$$\theta\left|_U\right.= \xi_1\frac{d x_1}{x_1}+\sum_{i=2}^n\xi_idx_i$$
and $\gamma^{-1}(Z)\cap U=\{x_1=0\}.$
There is a holomorphic map 
$\delta :\{ t\in \mathbb{C}: |t|<1 \}\to\gamma^{-1}(\Gamma)$
such that 
$$\gamma(\delta
(0))=o \text{ and } \delta^{-1}(\gamma^{-1}( \pi^{-1}(N)  ))=\{0\}.$$
Set
$\delta_i=x_i\circ\delta$, $1\le i\le n$. Since $\theta$ vanishes
on $\gamma^{-1}(\Gamma_0\setminus  Z)$,

\begin{equation}\label{argument}
\xi_1(\delta(t))\frac{\delta_1'(t)}{\delta_1(t)}+
\sum_{i=2}^n\xi_i(\delta(t))\delta_i'(t)=0 \qquad \hbox{\rm if} ~
t\not=0.
\end{equation}
Hence 
\begin{equation}\label{CSI}
\xi_1(o)=0,
\end{equation}
and (\ref{contido}) holds.

\noindent Since $Z$ is invariant, $\Gamma\cap Z$ is an involutive
submanifold of $P^*\langle M/N\rangle$. Hence $\Gamma\cap Z$ is an
involutive submanifold of $\mathbb{P}^*\langle Q/R\rangle$. Hence its
irreducible components are involutive. Since dim$\Gamma_0=$
dim$\Gamma-1$, $\Gamma_0\setminus \pi^{-1}(R)$ is a Legendrian
subvariety of $\mathbb{P}^*( Q\setminus R)$. Let $S_0$ be the closure
in $Q$ of the projection of $\Gamma_0\setminus \pi^{-1}(R)$.
Then $\Gamma_0$ is the conormal of $S_0$. By the induction
hypothesis, $\Gamma_0$ is a Legendrian variety of $\mathbb{P}^*\langle
Q/R\rangle$.
\end{proof}

\begin{thm}\label{colombo}
An irreducible Legendrian subvariety of a projective logarithmic
cotangent bundle is the conormal of its projection.
\end{thm}
\begin{proof}
The result is known for Legendrian subvarieties of a projective
cotangent bundle (see for instance \cite{PHAM}). The theorem is an
immediate consequence of this particular case.
\end{proof}

\section{Blow ups}

\noindent
Let $D$ be a submanifold of $M$. 
The vector bundle $T_DM$ defined by the exact sequence of vector bundles
$
0\rightarrow TD \rightarrow D\times_M TM \rightarrow T_DM \rightarrow 0
$
is called the  \emph{normal bundle} of $M$ along $D$.

\begin{lemma}\label{canonicaltangent}
Let $f:X \rightarrow Y$ be a holomorphic map between  manifolds. Let $A$ $[B]$ be a submanifold of $X$ $[Y]$. If $f(A)=B$ and $f$ and $f\left|_A \right.:A\rightarrow B$ are submersions, $f$ induces a canonical holomorphic map $\sigma$ from $T_AX$ into $T_BY$.
\end{lemma}
\begin{proof}
Given $o\in X$, $Df(a)$ defines maps from $T_oX$ onto $T_{f(o)}Y$ and from $T_oA$ onto  $T_{f(o)}B$. Hence $Df(o)$ induces a map from $T_oX/T_oA$ onto \linebreak $T_{f(o)}X/T_{f(o)}B$. Therefore $Df$ induces a map $\sigma:T_AX \rightarrow T_BY$. Locally there are coordinates 
$(x_1,\ldots ,x_a,y_1,\ldots,y_b,z_1,\ldots,z_c,w_1,\ldots,w_d)$
on $X$ and  $(u_1,\ldots ,u_a,v_1,\ldots,v_c)$ on $Y$ such that
$A=\{z=w=0\}, B=\{v=0\}$
and $f(x,y,z,w)=(x,z).$ 
Hence there are local coordinates
$$(x_1,\ldots ,x_a,y_1,\ldots,y_b,\widetilde{z}_1,\ldots,\widetilde{z}_c,\widetilde{w}_1,\ldots,\widetilde{w}_d)$$
on $T_AX$ and $(u_1,\ldots ,u_a,\widetilde{v}_1,\ldots,\widetilde{v}_c)$ on $T_BY$ 
such that $A=\{\widetilde{z}=\widetilde{w}=0\}$, $B=\{\widetilde{v}=0\}$ and 
$\sigma(x,y,\widetilde{z},\widetilde{w})=(x,\widetilde{z}).$
\end{proof}

\noindent
Let $X^\#_M$ be the normal deformation of $M$ in $X$ (see section 4.1 of \cite{ka3}).
Let $\tau:\widetilde X_M\to X$ be the blow up of $X$ along $M$.
Set $E=\tau^{-1}(M)$.
There are maps $\Phi:X^\#_M \to \widetilde X_M$, $p:X^\#_M \to X$ and $s:X^\#_M \to \mathbb C$ such that:

\noindent
$(i)$ $p^{-1}(X\setminus M)\simeq (X\setminus M)\times \mathbb C^*$;

\noindent
$(ii)$ $s^{-1}\mathbb C^*\simeq X\times \mathbb C^*$, $s^{-1}(\{0\})\simeq T_MX$;

\noindent
$(iii)$ $p=\tau\Phi$.

\noindent
There is a free action of $\mathbb C^*$ on $X^\#_M$ such that $\Phi$ induces an isomorphism from 
$X^\#_M/\mathbb C^*$ into $\widetilde X_M$ that takes $\mathbb P_MX$ into $E$.

\noindent
Let $\widetilde S$ [$\widetilde C_D(S)$] be the proper inverse image of $S$ by $\tau$ [$p$].
The set $C_D(S)=\widetilde C_D(S)\cap T_MX$ is called the \em normal cone of $S$ along $M$. \em
We recall that $\Phi (C_D(S))=\widetilde S\cap E$.

\noindent 
Set $X=\mathbb{C}^{a+b+c}$ with coordinates
$(x_1, \ldots, x_a, y_1, \ldots, y_b, z_1, \ldots, z_c).$

\noindent 
Set $\Lambda=\{x=y=0\}$ and set
$L=\{(\tilde{x},\tilde{y},z)\in T_{\Lambda}X: \widetilde{x}=0\}.$
The blow up of  $X$ along $\Lambda$ is the glueing of the affine open sets $U_{x_i}, 1 \leq i \leq a, U_{y_j}, 1 \leq j \leq b.$
Here $U_{x_i}$ is the affine set with coordinates 
$(\frac{x_1}{x_i},\ldots,\frac{x_{i-1}}{x_i},x_i,\frac{x_{i+1}}{x_i},$ $\ldots,\frac{x_{a}}{x_i},
\frac{y_{1}}{x_i},\ldots,\frac{y_b}{x_i},z_1,\ldots,z_c)$.

\begin{lemma}\label{key}
Let $\Gamma$ be the germ of a closed analytic subset of $X$. 

\noindent
If $C_{\Lambda}(\Gamma)\cap L \subset \{\tilde{x}=\tilde{y}=0\}$,
$\widetilde{\Gamma}\cap E\subset \cup_{i=1}^aU_{x_i}.$
\end{lemma}
\begin{proof} 
Notice that ${\mathbb C^{a+b+c}_\Lambda}^{\#}=\{(s,\widetilde x,\widetilde y,z):(\widetilde x,\widetilde y)\neq (0,0)\}$. Moreover, 
$\Phi$ induces a surjective map
\[
(s,\widetilde x, \widetilde y,z)\mapsto 
\left(\frac{\widetilde x_1}{\widetilde x_i},\ldots,s\widetilde x_i,\ldots,\frac{\widetilde x_a}{\widetilde x_i},
\frac{\widetilde y_1}{\widetilde x_i},\ldots,\frac{\widetilde y_b}{\widetilde x_i},z_1,\ldots,z_c \right)
\]
from $\{ (s,\widetilde x, \widetilde y,z): \widetilde x_i\neq 0 \}$ into $U_{x_i}$.
Hence $E\cap U_{x_i}=\Phi(\{ s=0, \widetilde x_i\neq 0 \})$ and
$E\setminus U_{x_i}=\Phi(\{ s=0, \widetilde x_i= 0 \})$.
Since $E\setminus \cup_{i=1}^aU_{x_i}=\Phi(\{ s=0, \widetilde x= 0 \})$,
\[
(\widetilde \Gamma \cap E)\setminus \cup_{i=1}^aU_{x_i}=\Phi(C_\Lambda(\Gamma) \cap L ).
\]
\end{proof}

\noindent Let $L$ be a submanifold of a manifold $M$ of codimension greater or equal than $2$. Let $N$ be a normal crossings divisor of $M$. We say that $L$ is a \em center \em of $(M,N)$ at $o \in L \cap N$ if there are manifolds $\Sigma_1, \ldots, \Sigma_s$ and $I,J \subset \{1, \ldots, s\}$ such that  $\Sigma_1\cup \ldots \cup \Sigma_s$ is a divisor with normal crossings at $o$ and 
$$(N,o)=(\cup_{i \in I}\Sigma_i,o), (L,o)=(\cap_{i \in J}\Sigma_i,o).$$
We say that $L$ is a \em trivial center \em [\em non-trivial center\em] of $(M,N)$ at $o$, if $J \subset I [J \not \subset I]$.

\begin{thm}\label{inversecommutes}
Let $M$ be a  manifold and let $N$ be a divisor with normal
crossings of $M$. Let $L$ be a trivial center of $(M,N)$.
Let $\rho:\widetilde M\to M$ be the blow
up of $M$ along $L$. Set $\widetilde N=\rho^{-1}(N)$. 

\noindent \em(i)\em The
blow up of $T^*\langle M/N\rangle$ along $\pi^{-1}(L)$ is a
logarithmic symplectic manifold isomorphic to $T^*\langle \widetilde
M/\widetilde N \rangle$ and diagram \em (\ref{dia1}) \em commutes.
\begin{equation}\label{dia1}
\begin{array}{ccccc}
  T^*\langle M/N\rangle  & \leftarrow &  T^*\langle \widetilde
M/\widetilde N \rangle  \\
  \downarrow &  &   \downarrow \\
  M  & \leftarrow &   \widetilde M
\end{array}
\end{equation}

\noindent 
\noindent \em(ii)\em If $S$ is a
natural hypersurface of $M$, the proper inverse image of the
conormal of $S$ equals the conormal of the proper inverse image of
$S$.
\end{thm}

\begin{proof}
The proof of statement (i) is similar to the proof of statement (ii) of Theorem \ref{diagteo}. Hence we omit it. Assume $S$ irreducible. Set $\Gamma=T_S^*\langle M\setminus N \rangle$. Let $\widetilde{\Gamma}$ be the proper inverse image of $\Gamma$ by the blow up with center $\pi^{-1}(L)$. Since diagram (\ref{dia1}) commutes, the projection of $\widetilde{\Gamma}$ into $\widetilde{M}$ equals the proper inverse image $\widetilde{S}$ of $S$ by $\rho$. Since $\widetilde{\Gamma}\setminus \pi^{-1}(\widetilde{N})$ is the conormal of  $\widetilde{S}\setminus \widetilde{N}$, $\widetilde{\Gamma}=T_{\widetilde{S}}\langle \widetilde{M} \setminus \widetilde{N} \rangle.$ \end{proof}

\noindent Let $X$ be a  manifold and let $Y$ be a closed hypersurface of
$X$. We will denote by $\mathcal O_{X}\gr{Y}$ the sheaf of
meromorphic functions $f$ such that $fI_Y\subset\mathcal O_{X}$.

\begin{theorem}\label{diagteo}
Let $N$ be the normal crossings divisor of a  manifold $M$. Let $L$ be a nontrivial center of $(M,N)$. Let $\tau$ be the blow up of $X=\mathbb{P}^*\langle M/N \rangle$ along $\Lambda=\mathbb{P}_L^*\langle M/N \rangle$. Set $E=\tau^{-1}(\Lambda)$. Let $\rho:\widetilde{M}\rightarrow M$ be the blow up of $M$ along $L$. Set $\widetilde{N}=\rho^{-1}(N)$.

\noindent
\em (i) \em If $\mathcal{L}$ is the canonical contact structure of $\mathbb{P}^*\langle M/N \rangle$, the $\mathcal{O}_{\widetilde{X}}$-module $\mathcal{O}_{\widetilde{X}}(E)\tau^*\mathcal{L}$ is a structure of logarithmic contact manifold on $\widetilde{X}$ with poles along $\tau^{-1}(\pi^{-1}(M))$.

\noindent \em (ii) \em There is an injective
contact transformation $\varphi$ from a dense open subset $\Omega$ of $\widetilde X$ onto $P^*\langle \widetilde M/\widetilde N \rangle$ such that
diagram \em (\ref{dia2}) \em commutes.

\begin{equation}\label{dia2}
\begin{array}{ccccc}
  P^*\langle M/N\rangle  & \overset{\tau}{\leftarrow} &
  \widetilde X \hookleftarrow  \Omega &
  \overset{\varphi}{\hookrightarrow} & P^*\langle \widetilde M/\widetilde N \rangle  \\
  \!\!\!\!\!\! \pi \downarrow & &  & &   \downarrow \pi \\
  M  & &  \overset{\rho}{ \longleftarrow}   &  & \!\!\! \widetilde M
\end{array}
\end{equation}

\noindent
\em (iii) \em Let $S$ be a germ of a natural analytic subset of $(M,N)$ at $o \in N$. Set $\Gamma=\mathbb{P}^*_S\langle M/N\rangle$. Let $\widetilde{S}$ be the proper inverse image of the blow up of $M$ along $L$. If $S$ has trivial limits of tangents at $o$ and $C_{\Lambda}(\Gamma)\cap\sigma^{-1}(L)\subset\Lambda$, then $\widetilde{\Gamma} \subset \Omega$ and $\varphi(\widetilde{\Gamma})=\mathbb{P}^*_{\widetilde{S}}\langle \widetilde{M}/\widetilde{N}\rangle$, where $\sigma$ denotes the canonical projection from $T_{\Lambda}\mathbb{P}^*\langle M/N\rangle$ onto $T_LM$ introduced in Lemma \ref{canonicaltangent}.
\end{theorem}

\proof
The proof of statement $(i)$ is quite similar to the proof of Proposition 9.2 of \cite{neto4}. Hence we omit it.

\noindent 
(ii) Assume $M= \mathbb{C}^{n+1}$, $N=\{x_1\cdots x_{\nu}=0\}$ and $X=\mathbb{P}^*\langle M/N \rangle$. 

\noindent The canonical 1-form $\theta$ of $T^*\langle M/N \rangle$ equals 
$$\sum_{i=1}^{\nu} \xi_i \frac{dx_i}{x_i}+\sum_{i=\nu+1}^{n+1}\xi_i dx_i.$$  
Hence there is $\iota$ such that $L=\{x_{\iota}=\cdots=x_k=x_{n+1}=0\}$.
\noindent Let $\widehat{\widetilde{X}}$ be the homogeneous symplectic manifold associated to $\widetilde{X}$. Let $\widehat{\theta}$ be the canonical 1-form of $\widehat{\widetilde{X}}$. By the argument of (i) $\widehat{\widetilde{X}}$ is the union of open set $\widehat{U}_j, j=\iota,\cdots, k,n+1$ and $\widehat{V}_j, j=1,\ldots, \nu,k+1, \cdots, n$. 

\noindent Set $\widehat{\Omega}=\cup_j \widehat{U}_j $. Set $\Omega=\widehat{\Omega}_*$, $\widehat{\theta}_j=\widehat{\theta} \left. \right |_{\widehat{U}_j}$, $j=2, \ldots, k, n+1$.  Endow $\mathbb{C}^{2n+2}$ with the coordinates
$$x_1, \ldots, x_{\iota-1},{\textstyle\frac{x_{\iota}}{x_j}},\ldots,{\textstyle\frac{x_{j-1}}{x_j}},x_j,{\textstyle\frac{x_{j+1}}{x_j}},\ldots, {\textstyle\frac{x_{\nu}}{x_j}},x_{\nu+1},\ldots,x_{n}, {\textstyle\frac{x_{n+1}}{x_j}}       ,\eta_1, \ldots, \eta_{n+1}.$$
We can assume that  $\widehat{U}_j=\{(\eta_1,\ldots,\eta_{n+1})\neq(0,\cdots,0)\}$ and
\begin{eqnarray*}
\widehat{\theta}_j=\sum_{i=1}^{\iota-1}\eta_i\frac{dx_i}{x_i}+\sum_{i=\iota}^{\nu}\eta_i\frac{d\frac{x_i}{x_j}}{\frac{x_i}{x_j}}+\eta_jdx_j+\sum_{i=\nu+1}^k\eta_id{\textstyle\frac{x_i}{x_j}}+\sum_{i=k+1}^{n}\eta_idx_i+\eta_{n+1}d{\textstyle\frac{x_{n+1}}{x_j}}.
\end{eqnarray*}
The blow up of $M$ along $L$ is the glueing of the open affine sets $M_j, j=\iota, \ldots,k,n+1$ where $M_j$ is associated to the generator $x_j$ of the defining ideal $L$. Hence $T^*\langle \widetilde{M}/\widetilde{N} \rangle$ is the glueing of the open affine sets $T^*\langle M_j/\widetilde{N}\cap M_j\rangle, j=\iota,\ldots,k,n+1.$

\noindent Set $\widehat{W}_j=\ballontop{T}^*\langle M_j/N\cap M_j \rangle$. Let $\widetilde{\theta}$ be the canonical 1-form of $T^*\langle M/N \rangle$. Set $\widetilde{\theta}_j=\left. \widetilde{\theta} \right|_{\widehat{W}_j}$.
Endow $\mathbb{C}^{2n+2}$ with the coordinates
$$x_1, \ldots, x_{\iota-1},{\textstyle\frac{x_{\iota}}{x_j}},\ldots,{\textstyle\frac{x_{j-1}}{x_j}},x_j,{\textstyle\frac{x_{j+1}}{x_j}},\ldots,{\textstyle\frac{x_{\nu}}{x_j}},x_{\nu+1},\ldots,x_{n}, {\textstyle\frac{x_{n+1}}{x_j}} ,\zeta_1, \ldots, \zeta_{n+1}.$$

\noindent We can assume that 
$\widehat{W}_j=\{(\zeta_1,\ldots,\zeta_{n+1})\neq(0,\cdots,0)\}$
and
\begin{eqnarray*}
\widetilde{\theta_j}=\sum_{i=1}^{\iota-1}\zeta_i\frac{dx_i}{x_i}+\sum_{\substack {i=\iota \\ i \neq j}}^{\nu}\zeta_i\frac{d\frac{x_i}{x_j}}{\frac{x_i}{x_j}}+\zeta_j {\textstyle \frac{dx_j}{x_j}}+\sum_{\substack{i=\nu+1 \\ i \neq j}}^k\zeta_id\frac{x_i}{x_j}+\sum_{i=k+1}^{n}\zeta_idx_i+\zeta_{n+1} d{\textstyle\frac{x_{n+1}}{x_j}}.
\end{eqnarray*}
Since 
$$\widehat{\widetilde{X}} \hookleftarrow T^*\langle M \setminus L / N \setminus L\rangle \simeq T^*\langle \widetilde{M} \setminus \rho^{-1}(L) / \widetilde{N} \setminus \rho^{-1}(L)\rangle \hookrightarrow T^*\langle \widetilde{M}\setminus \widetilde{N}\rangle$$
There is a bimeromorphic contact transformation $\widehat{\varphi}^{-1}:\widehat{\widetilde{X}} \rightarrow \ballontop{T}^*\langle \widetilde{M}/\widetilde{N} \rangle.$
It is enough to show that the domain of $\widehat{\varphi}$ contains $\widehat{\Omega}$ and its image equals $\ballontop{T}^*\langle \widetilde{M}/\widetilde{N} \rangle$. Since 
$$\widehat{U}_j\setminus \tau^{-1}(\pi^{-1}(L))=\widehat{W}_j\setminus \pi^{-1}(\rho^{-1}(L)),$$
$\eta_i=\zeta_i$ on a dense open set of their domain. Hence $\eta_i=\zeta_i$ everywhere and the domain of $\widehat{\varphi}$ contains $U_j$ for $j=\iota, \ldots, k,n+1$.

\noindent (iii) The statement follows from Lemma \ref{key} and the arguments of the proof of statement (ii) of theorem \ref{inversecommutes}.  
\qed

\newpage
\section{Main Result}
\noindent 
Let $N$ be a germ of normal crossings divisor of a three dimensional manifold $M$ at a point $o$. 
Let $S$ be a germ of a natural irreducible  surface of $M$ at $o$. 
We say that $S$ is \em quasi-ordinary \em at $o$ if there is a system of local coordinates 
$(x,y,z)$ on a neighbourhood of $o$ such that the discriminant of $S$ relatively to the projection 
$(x,y,z)\mapsto (x,y)$ is contained in $\{xy=0\}$. 
There is a positive integer $m$ and 
$\zeta \in \mathbb C \{x,y\}$ such that 
$z=\zeta ( x^{\frac{1}{m}},y^{\frac{1}{m}} )$ 
defines a parametrization of $S$. 
If $\zeta \neq 0$, let $x^{\lambda}y^{\mu}$ be the monomial of lowest degree that occurs in $\zeta$. 
If $\lambda\mu=0$ or $\lambda,\mu \in \mathbb{Z}$ and there is a monomial of $\zeta$ that does not verify this condition, 
let $x^ay^b$ be the monomial of lowest degree of $\zeta$ that does not verify this condition. Otherwise, set $a=b=+\infty$. 
If $\zeta=0$, we set also $\lambda=\mu=+\infty$. We call $\lambda, \mu, a, b$, the \em exponents \em of $\zeta$.

\noindent 
We say that $\zeta$ is in \em strong normal form \em if $\zeta$ has no monomials with integer exponents and 
\begin{equation}\label{cond}
\lambda>\mu ~ \hbox{ or } ~ \lambda=\mu \hbox{ and } a\geq b. 
\end{equation}
We say that $\zeta$ is in \em normal form \em if $\zeta=0$ or (\ref{cond}) holds.

\noindent
Lipman proved a desingularization theorem for quasi-ordinary surfaces (see \cite{Lipmanphd}).
It was not clear that Lipman's result produced an embedded desingularization 
because the parametrization is related to a system of local coordinates and the choice of the centers
is dependent on another system of local coordinates, related to the divisors created by successive explosions.
Ban and Mcewan showed that Lipman's algorithm produces an embedded desingularization that coincides with the
general 
algorithm proposed by Bierstone and Milman \cite{milman}.
Theorem 2.5 of \cite{mcewan} shows that after all the two systems of local coordinates are not that different.
This result is a key point of their proof. Theorem \ref{coord} is a slightly more precise version of theorem 2.5 of \cite{mcewan}.
We will need it in order to prove Theorem \ref{loglimitsteo}.

\begin{definition} 
Let $M$ be a  manifold of dimension $3$.
Let $N$ be a normal crossings divisor of $M$.
Let $S$ be a natural surface of $(M,N)$.
Let $o\in S\cap N$. We say that $N$ is \em adapted to \em $S$ \em at \em $o$ if there is
a system of local coordinates $(x,y,z)$ centered at $o$ such that $(N,o)\subset \{xyz=0\}$, the discriminant of the germ of $S$ relatively to the projection $(x,y,z)\mapsto (x,y)$ is contained in $\{ xy=0\}$ and the parametrization $z=\varphi (x^{\frac{1}{m}},y^{\frac{1}{m}})$ of $(S,o)$ is in normal form. We say that $N$ is \em adapted to \em $S$ if $N$ is adapted to $S$ at $o$, for each $o\in S\cap N$.
\end{definition}

\noindent
We will denote by $\varepsilon$ or $\varepsilon_i$ a unit of $\mathbb{C}\{x^{\frac{1}{m}},y^{\frac{1}{m}}\}$ and by $\delta$ or $\delta_i$ a unit of $\mathbb{C}\{x^{\frac{1}{m}}\}$, for a convenient $m$.

\begin{theorem}\label{coord}
Let $S_0$ be the germ of a quasi-ordinary surface.
At each step of 
Lipman's algorithm the normal crossings divisor
we obtain is adapted to the proper inverse
image of $S$.
\end{theorem}

\begin{proof} 
We prove the theorem by induction on the number of steps. Let $S$ be the proper inverse image of $S_0$ at some step of Lipman's algorithm and $N$ the system of exceptional divisor at that step. Let $o\in S\cap N$. We assume that $S$ is singular at $o$. Let $(x,y,z)$ a system of coordinates centered at $o$ such that $S$ is defined at $o$ by the parametrization $z=x^{\lambda}y^{\mu}\varepsilon$ and $(N,o)\subset \{xyz=0\}$. For this proof we will only consider the following types of changes of coordinates:
\begin{eqnarray}
(x,y,z)\mapsto(y,x,z),\qquad\qquad\qquad\quad\;\label{T1}\\
(x,y,z)\mapsto(z,y,x),\qquad\qquad\qquad\quad\;\label{T2}\\
(x,y,z)\mapsto(z,x,y),\qquad\qquad\qquad\quad\;\label{T3}\\
(x,y,z)\mapsto(x,y,z-q),\;q\in\mathbb C\{x,y\}.\label{T4}
\end{eqnarray}
We will use the notation of table \ref{genposi} for the center of the resolution step.

\noindent
Assume $\lambda+\mu<1$. Then we blow up $\sigma_0$. Let $\widetilde S$, $\widetilde{N}$ be the proper inverse images of, respectively, $S$ and $N$. Set $U_x$ as the open affine chart defined by the coordinates $(x,\frac{y}{x},\frac{z}{x})$. Then, after a re-parametrization,
\begin{equation}\label{ADAPTED1}
\frac{y}{x}=\left(\frac{z}{x}\right)^{\frac{1}{\mu}}x^{\frac{1-(\lambda+\mu)}{\mu}}\varepsilon_1
\end{equation}
and $\widetilde{N}\cap U_x\subset \{x\frac{y}{x}\frac{z}{x}=0\}$. Notice that $\frac{1}{\mu}>\frac{1-(\lambda+\mu)}{\mu}$ and $\frac{1}{\mu}>1$. Furthermore, if $(\lambda,\mu)=(a/n,1/n)$, for some positive integers $a,n$, $\frac{1}{\mu}$ and $\frac{1-(\lambda+\mu)}{\mu}$ are positive integers. After a change of coordinates of type (\ref{T3}), $\widetilde{N}$ is adapted to $\widetilde{S}$ at the origin of $U_x$. From (\ref{ADAPTED1}), one can easily check that, after a appropriate change of coordinates, $\widetilde{N}$ is adapted to $\widetilde{S}$ at the points of $U_x$ of the type $(0,0,a)$ and $(a,0,0)$, $a\in\mathbb C^{\ast}$. On the open affine chart defined by the coordinates $(\frac{x}{y},y,\frac{z}{y})$ the reasoning is analogous. On the open affine chart defined by the coordinates $(\frac{x}{z},\frac{y}{z},z)$, no re-parametrization is needed and one can easily check that $\widetilde{N}$ is adapted to $\widetilde{N}$ on this chart. Notice that if on this chart the multiplicity drops, we return to the case $\lambda+\mu\geq 1$.

\noindent
Assume $\lambda+\mu\geq 1$. Let $p\in\mathbb C\{x,y\}$ such that $z-p$ has no monomials with integer exponents. Assume that for $z-p$, $1<\lambda<2$ and $\mu=0$. Then we blow up $\sigma_x$. On the open affine set $V_x$  defined by the coordinates $(x,y,\frac{z}{x})$, possibly after a change of coordinates of the type (\ref{T4}) with $q\in\mathbb C^{\ast}$, we need to do a re-parametrization and we obtain
\begin{displaymath}
x=\left(\frac{z}{x}\right)^{\frac{1}{1-\lambda}}\varepsilon_1
\end{displaymath}
and $\widetilde{N}\cap U_x\subset \{xy\frac{z}{x}=0\}$. Notice that if $\lambda=\frac{n-1}{n}$, for some positive integer $n$, $\frac{1}{1-\lambda}$ is a positive integer. From this parametrization one can easily check that $\widetilde{N}$ is adapted to $\widetilde{S}$. The remaining cases need no re-parametrization and one can easily check that the theorem holds.

\noindent
We remark that if $\lambda,\mu\in\mathbb Z$, at some step of Lipman's algorithm, before we reach the situation where $\lambda+\mu<1$, we do succession of resolution steps with center $\sigma_x$ or $\sigma_y$. Eventually we reach a situation where we can do a change of coordinates of type (\ref{T4}) with $q=a+p$, $a\in\mathbb C^{\ast}$. Hence, for $o\in S$ and an adequate system of local coordinates $(x,y,z)$ centered at $o$, $(N,o)\subset\{xy=0\}$ and we can consider a parametrization of $(S,o)$ that is in strong normal form.
\end{proof}

\begin{theorem}\label{loglimitsteo}
Let $M$ be a  manifold of dimension $3$.
Let $S$ be a quasi-ordinary surface of $M$ at $o$.
Let $N$ be a general configuration adapted to $S$ at $o$.
Let $\Sigma$ be the logarithmic limit of tangents of $S$ relatively to $N$. 

\noindent
$(a)$
Assuming that $p+\zeta$ is in normal form:

\noindent \em (i) \em If $T_oN=\emptyset$, $\Sigma$ is trivial if and only if 
   
\em ($\phi$1) \em $\mu\geq 1$ or

\em ($\phi$2) \em  $\mu=0$, $b \geq 1$.

\noindent \em (ii) \em If $T_oN=T_o\{x=0\}$, $\Sigma$ is trivial if and only if

\em (x1) \em $\mu\geq 1$ or

\em (x2) \em  $\mu=0$ and $b \geq 1$.

\noindent \em (iii) \em If $T_oN=T_o\{y=0\}$, $\Sigma$ is trivial if and only if

\em (y1) \em $\lambda\geq 1$.

\noindent \em(iv) \em If $T_oN=T_o\{z=0\}$, $\Sigma$ is trivial if and only if

\em (z1) \em $\mu=0$, $b\geq 1$. 

\noindent \em(v) \em If  $T_oN=T_o\{xz=0\}$, $\Sigma$ is trivial if and only if

\em (xz1) \em $ \mu \neq 0$ or

\em (xz2) \em  $\mu=0$, $b \geq 1$ or

\em (xz3) \em  $\mu=0$, $b<1$, $a=\lambda$.

\noindent \em(vi) \em $\Sigma$ is trivial if $T_oN=T_o\{xy=0\}$, $T_oN=T_o\{yz=0\}$ or $T_oN=T_o\{xyz=0\}$.
  
\vskip .1in
  
\noindent
$(b)$ Let $S$ be the germ of a quasi-ordinary surface at $o$ with trivial limit of tangents. 
Set $\Gamma=\mathbb{P}^*_S\langle M/N \rangle$. Let $L$ be one of the admissible centers for $S$ considered in Table \em\ref{genposi}. \em 
Set $\Lambda=\mathbb{P}^*_L\langle M/N \rangle$. 
Let $\widetilde{\Gamma}$ be the proper inverse image of $\Gamma$ by the blow up of $\mathbb{P}^*\langle M/N\rangle$ with center $\Lambda$. 
Then 
$$
\widetilde{\Gamma}=\mathbb{P}^*_{\widetilde{S}}\langle M/N \rangle.
$$
\end{theorem}

\begin{proof}
Let $S$ be a surface that verifies the conditions of statement (b).
Let $\lambda$ be the only limit of tangents of $S$ at $o$. 
Let $\sigma:T_{\Lambda}\mathbb{P}^*\langle M/N \rangle \rightarrow T_LM$
be the map associated to $\pi:\mathbb P\langle M/N\rangle\to M$ by Lemma \ref{canonicaltangent}.
By Lemma \ref{key} it is enough to verify that  
\begin{equation}\label{ALG}
C_{\Lambda}(\Gamma)\cap\sigma^{-1}(L)\subset \Lambda
\end{equation}
holds in order to show that $\widetilde{\Gamma}\subset \Omega$, where $\Omega$ is the set referred to in diagram $\ref{dia2}$.
We will prove statements (a) and (b) in each of the cases considered in statement (a). The cases where $\zeta$ does not depend on $y$ are easy to handle so we omit them.

\noindent 
In case (i) the triviality of the limits of tangents  was already discussed in \cite{LIMITS}. Set $\theta=\xi dx+\eta dy+\zeta dz=\zeta(dz-pdx-qdy)$. 
In case ($\phi 1$), $L=\{x=y=z=0\}$. Moreover, the blow up of $\mathbb P^*M$ along $\Lambda=\{x=y=z=0\}$ equals $\mathbb P^*\langle \widetilde M/E\rangle$.
Hence (\ref{ALG}) is trivially verified.

\noindent
In case ($\phi 2$) $L=\{x=z=0\}$. Hence $\Lambda=\{x=z=q=0\}$ and $\sigma(\widetilde{x},y,\widetilde{z},p,\widetilde{q})=(\widetilde{x},y,\widetilde{z}).$
Since $\mu=0$,
\begin{equation}\label{v2_1}
a \geq \lambda \geq 1.
\end{equation}
 Since 
$z=a_{\lambda 0}x^{\lambda}+\ldots +a_{a b}x^{a}y^{b}+\cdots$,
\begin{equation}\label{v2_2}
q=\frac{\partial z}{\partial y}=x^{a}y^{b-1}\varepsilon
\end{equation}
It follows from $(\ref{v2_1})$ and $(\ref{v2_2})$ that $\Gamma$ is contained in a hypersurface 
$q^n+\sum_{i=0}^{n-1}a_iq^i=0$
where $a_i \in \mathbb{C}\{x,y\}$ and $a_i \in (x)^{n-i}$.
Hence there are $\widetilde{a}_i \in (\widetilde{x})$ such that $C_{\Lambda}(\Gamma)$ is contained in an hypersurface 
$\widetilde{q}^n+\sum_{i=0}^{n-1}\widetilde{a}_i\widetilde{q}^i=0$.
Therefore (\ref{ALG}) holds.

\noindent 
In case (ii) $\theta=\zeta (dz-p{dx}/{x}-qdy)=\xi dx/x+\eta dy + \zeta dz.$

\noindent 
Assume  $0 <\mu<1$. 
Then $z=x^{\lambda}y^{\mu}\varepsilon_1$ is a parametrization of $S$ and 
$$z=x^{\lambda}y^{\mu}\varepsilon_1, ~~~ p=x^{\lambda}y^{\mu}\varepsilon_2, ~~~ q=x^{\lambda}y^{\mu-1}\varepsilon_3.$$
is a parametrization of the regular part of $\Gamma_{reg}$.

\noindent Set 
$\beta={\lambda \alpha}/({1-\mu}),$
where $\alpha$ is a positive integer.
\noindent There are $A, B \in \mathbb{C}^*$, and units $\sigma_i$ of $\mathbb{C}\{t\}$, $1 \leq i \leq 3$, such that the map that takes $t$ into
$$
\left(At^{\alpha}, Bt^{\beta},A^{\lambda}B^{\mu}t^{\alpha\lambda+\beta\mu}\sigma_1,A^{\lambda}B^{\mu}t^{\alpha\lambda+\beta\mu}\sigma_2, A^{\lambda}B^{\mu-1}t^{\alpha\lambda+\beta(\mu-1)}\sigma_3 \right)
$$
is a curve of $\Gamma$. Since $\alpha\lambda+\beta\mu>0$ and $\alpha\lambda+\beta(\mu-1)=0$, 
the curve converges to $(0,0,0,0:A^{\lambda}B^{\mu-1}\sigma_3(0):1)$ when $t$ goes to 0. Hence $\Sigma$ is not trivial.

\noindent 
Assume $\mu=0$ and $b<1$. Then
$z=x^{\lambda}\delta_1+x^{a}y^{b}\varepsilon_1$, $p=x^{\lambda}\varepsilon_2$, $q=x^{a}y^{b-1}\varepsilon_3$
define a parametrization of $\Gamma_{reg}$. Hence we can repeat the previous argument.

\noindent 
Assume $\mu \geq 1$. Hence $z=x^{\lambda}y^{\mu}\varepsilon_1$ defines a parametrization of $S$ and 
$ q=x^{\lambda}y^{\mu-1}\varepsilon_2$
defines a parametrization of a hypersurface that contains $\Gamma.$
Hence $\Sigma\subset \{\eta=0\}$. By (\ref{CSI}), $\Sigma \subset \{\xi=0\}$.

\noindent 
If $\mu=0$ and $b\geq 1$ we can obtain a proof of the triviality of $\Sigma$ combining the arguments of the previous cases.

\noindent 
In case (x1), $\mu \geq 1$. If $L=\{x=y=z=0\}$, $\Lambda=\{x=y=z=p=0\}$ and 
$\sigma(\widetilde{x},\widetilde{y},\widetilde{z},\widetilde{p},q)=(\widetilde{x},\widetilde{y},\widetilde{z}).$
Then $z=x^{\lambda}y^{\mu}\varepsilon_1$, and $p=x^{\lambda}y^{\mu}\varepsilon_2$. Since $\lambda+\mu>2$,
$C_{\Lambda}(\Gamma)\subset\{\widetilde{p}=0\}$. Hence, (\ref{ALG}) holds.

\noindent  
If $L=\{y=z=0\}$, $\Lambda=\{y=z=p=0\}$. If $\mu>1$, the argument is similar to the previous one.  

\noindent 
If  $\mu=1$, then
$z=x^{\lambda}y \varepsilon_1$, $p=x^{\lambda}y \varepsilon_2$, and $C_{\Lambda}(\Gamma)\cap \sigma^{-1}(L)\subset\{\widetilde{p}=0\}$.

\noindent
(x2) Since $L=\{x=z=0\}$, $\Lambda=\{x=z=p=q=0\}$. 
Hence 
$$z=x^{\lambda} \varepsilon_1=x^{\lambda}\delta+x^{a}y^{b}\varepsilon_2 , ~~~ p=x^{\lambda}\varepsilon_3, ~~~ q=x^{a}y^{b-1}\varepsilon_4.$$
Since $\mu=0$, $\lambda>1$. Therefore $C_{\Lambda}(\Gamma)\subset\{\widetilde{p}=0\}.$

\noindent
Since $a\geq \lambda >1$, $b\geq1$. Therefore $C_{\Lambda}(\Gamma)\subset\{\widetilde{q}=0\}.$

\noindent 
(iii) 
By (\ref{CSI}), $\Sigma \subset \{\eta=0\}$.
If $\mu=0$ and $\lambda>1$, $\Sigma\subset\{\xi=0\}$ by the arguments of case (ii). The same arguments hold if $\lambda\geq 1$ and $\mu>0$.
If  $\lambda<1$, the argument of the first case considered in (ii) shows that $\Sigma$ is not trivial. 
The proof of (b) is similar to the the one presented in case (ii).
%
%
\noindent 
(iv) Set $\theta=\xi dx+\eta dy +\zeta {dz}/{z} =\zeta ({dz}/{z} -pdx-qdy)=\xi\left(dx-rdy-s{dz}/{z}\right).$

\noindent Assume that $\mu \neq 0$. Then
$$
z=x^{\lambda}y^{\mu}\varepsilon_1, ~~~ p=x^{\lambda-1}y^{\mu}\frac{\varepsilon_2}{z}, ~~~ q=x^{\lambda}y^{\mu-1}\frac{\varepsilon_3}{z}
$$
is a parametrization of $\Gamma_{reg}$. 
Hence there is a curve on $\Gamma$ of the type
$$
t \mapsto \left(At^{\alpha}, Bt^{\alpha},A^{\alpha}B^{\alpha}t^{\alpha(\lambda+\mu)}\sigma_1,\frac{\sigma_2}{At^{\alpha}},\frac{\sigma_3}{Bt^{\alpha}}  \right),
$$
where $\sigma_i, 1 \leq i \leq 3$ are units of $\mathbb{C}\{t\}$. Since
$\left( {\sigma_2}({At^{\alpha}})^{-1}:{\sigma_3}({Bt^{\alpha}})^{-1}:1 \right)=(B\sigma_1:A\sigma_2:ABt^{\alpha})$
converges to $(B\sigma_1(0):A \sigma_2(0):0)$, $\Sigma=\{\zeta=0\}.$

\noindent 
Assume that $\mu=0$ and $b<1$. 
Setting 
$\beta=\alpha({\alpha+1-\lambda})({1-b})^{-1}$,
we can show that $\Sigma=\{\zeta=0\}$.

\noindent 
Assume that $\mu=0$ and $b \geq 1$. 
Then
$z=x^{\lambda}\delta+x^ay^b\varepsilon_1=x^{\lambda}\varepsilon_2$, $p=\varepsilon_3/x$, $q=x^{a-\lambda}y^{b-1}\varepsilon_4$
define a parametrization of $\Gamma_{reg}$. 
Moreover, $\Gamma$ is contained in the hypersurfaces defined by the equations
\begin{equation}\label{eq_iv}
x \xi+\varepsilon_4 \zeta=0, ~~~ \eta+x^{a-\lambda}y^{b-1}\varepsilon_5\zeta=0.
\end{equation}
Hence $\Sigma=\{\eta=\zeta=0\}.$


\noindent Let us assume that $\mu=0$, $b\geq 1$, and prove (b). Set $L=\{x=z=0\}$. Hence $\Lambda=\{x=z=r=s=0\}$. 
From (\ref{eq_iv}), 
$\Gamma$ is contained in the hypersurfaces with parametrizations given by
$s=\varepsilon_4^{-1} x$, $r+x^{a-\lambda}y^{b-1}\varepsilon_5 s=0.$
Hence $C_{\Lambda}(\Gamma)\cap\{\widetilde{x}=\widetilde{z}=0\}\subset\{\widetilde{r}=\widetilde{s}=0\}.$


\noindent 
(v) Set 
$\theta=\zeta \left(dz/z -pdx/x-qdy\right)=\eta\left(dy-rdx/x-sdz/z\right).$
 
\noindent 
The case $\mu=0$ and $b\geq 1$ is quite similar to the one in (iv).

\noindent 
Assume that $\mu=0, b<1, a=\lambda$. 
Then
$z=x^{\lambda}\varepsilon_1, p=\varepsilon_2, q=y^{b-1}\varepsilon_3$
defines a parametrization of $\Gamma_{reg}.$
Hence $\Gamma$ is contained in the analytic set
$\xi+\varepsilon_2 \zeta=0, ~~~ y^{1-b}\eta+\varepsilon_3\zeta=0.$
Therefore $\Sigma=\{\xi=\zeta=0\}$.

\noindent Assume that $\mu=0, b<1, a>\lambda$. Setting
$\beta=\alpha(a-\lambda)/(1-b),$
it can be shown by the previous methods that there is a $u \in \mathbb{C}^*$ such that $\Sigma\supset\{(u:v:1):v \in \mathbb{C}^*\}.$

\noindent
We now consider statement (b) for the case $\mu=0$. We have that $L=\Lambda=\{x=z=0\}$. This situation is solved by theorem \ref{inversecommutes}.

\noindent 
Assume $\mu\neq 0$. 
%
\noindent
Set $L=\{x=y=z=0\}$. Hence $\Lambda=\{x=y=z=r=s=0\}.$ Then
$z=x^{\lambda}y^{\mu}\varepsilon_1$, $p=x^{\lambda}y^{\mu}z^{-1}\varepsilon_2=\varepsilon_3$,  $q=x^{\lambda}y^{\mu-1}z^{-1}\varepsilon_4=(\varepsilon_5 y)^{-1}$ is a parametrization of $\Gamma_{reg}$, and $\Gamma \subset \{\xi+\varepsilon_3 \zeta=y\varepsilon_5\eta+\zeta=0.\}$. Hence 
$r=-\varepsilon_3s, s=y\varepsilon_5$. Hence $C_{\Lambda}(\Gamma)\cap\rho^{-1}(L)\subset\{\widetilde{r}=\widetilde{s}=0\}.$

\noindent
If $L=\{ y=z=0\}$, then $\Lambda=\{y=z=r=s=0\}$, and this case is solved in a similar fashion.


\noindent
(vi) If $N=\{xy=0\}$, arguments previously used show that $\Sigma$ is trivial.



\noindent
Set $\theta=\xi dx/x+\eta dy/y+\zeta dz=\zeta\left(dz-pdx/x-qdy/y\right).$
Set $L=\{x=y=z=0\}$. Hence $\Lambda=\{x=y=z=p=q=0\}.$

\noindent Then
$z=x^{\lambda}y^{\mu}\varepsilon_1,  p=z^{\lambda}y^{\mu}\varepsilon_2=z\varepsilon_3, q=x^{\lambda}y^{\mu}\varepsilon_4=z \varepsilon_5$.
Hence $\Gamma$ is contained in the hypersurfaces  
$p^k+\sum_{i=0}^{k-1}a_ip^i=0$, $q^l+\sum_{i=0}^{l-1}b_iq^i=0,$
where $a_i \in (z^{k-i}), b_i \in (z^{l-i})$. Hence (\ref{ALG}) holds.

\noindent Assume that $\mu=0$, $L=\{x=z=0\}$ or $\mu \geq 1$, $L=\{y=z=0\}$. In both cases $C_{\Lambda}(\Gamma)\subset\{\widetilde{p}=\widetilde{q}=0\}$ by the standard arguments.


\noindent 
In the case $N=\{yz=0\}$, $\Sigma$ is always trivial by the arguments of case $N=\{xz=0\}$. 



\noindent  If $\mu \neq 0$, we are in the situation of (xz1). 

\noindent  Assume that $\mu=0$. Set $\theta=\xi dx+\eta dy/y+\zeta dz/z=\xi\left(dx-rdy/y-sdz/z\right).$

\noindent Then
$
z=x^{\lambda}\varepsilon_1=x^{\lambda}\delta+x^ay^b\varepsilon_2$, 
$p=x^{\lambda-1}z^{-1}\varepsilon_3=(x\varepsilon_4)^{-1}$, $q=x^ay^bz^{-1}\varepsilon_5=x^{a-\lambda}y^b\varepsilon_6$, define a parametrization of $\Gamma_{reg}$.
Therefore
$ \Gamma \subset \{x\varepsilon_4 \xi+ \zeta=\eta+x^{a-\lambda}y^b\varepsilon_6\zeta=0\}$. Hence $\Gamma$ is contained in the hypersurfaces
\begin{equation}\label{follows}
s=x\varepsilon_4, ~~~ r+x^{a-\lambda}y^b\varepsilon_6 s=0.
\end{equation}
\noindent 
It follows from (\ref{follows}) that $C_{\Lambda}(\Gamma)\cap\sigma^{-1}(L)\subset\{\widetilde{r}=\widetilde{s}=0\}$ if $L=\{x=y=z=0\}$ or $L=\{x=z=0\}$.


\noindent 
If $N=\{xyz=0\}$, $\theta=\xi {dx}/{x}+\eta {dy}/{y} +\zeta {dz}/{z} =\zeta ({dz}/{z} -p{dx}/{x}-q{dy}/{y})$.

\noindent 
Assume that $\mu \neq 0$.  $\Gamma$ is contained in the hypersurfaces determined by
$z=x^{\lambda}y^{\mu}\varepsilon_1$, $p=\lambda x^{\lambda}y^{\mu}\varepsilon_2z^{-1} \text{ and } q=\mu x^{\lambda}y^{\mu}\varepsilon_3z^{-1},$ 
where $\varepsilon_1(0)=\varepsilon_2(0)=\varepsilon_3(0)$. Hence 
$\Sigma=\{(\lambda:\mu:1)\}$.
A similar argument shows that we arrive to the same conclusion when $\mu=0$. 

\noindent In the case $N=\{xyz=0\}$, (b) is trivially verified.


\end{proof}

\begin{table}[!h]
\centering
\begin{tabular}{|l|l|l|l|}
\hline
Divisor & Conditions & Label & Center \\
\hline
\multirow{2}{*}{$\emptyset$}
& $\mu\geq 1$ & $\phi 1$ & $\sigma_0$ \\ 
& $\mu=0$ and $b\geq 1$ & $\phi 2$ & $\sigma_x$\\
\hline
\multirow{2}{*}{$\{x=0\}$}
& $\mu\geq 1$ & x1 & $\sigma_0$, $\sigma_y$\\ 
& $\mu=0$ and $b\geq 1$ & x2 & $\sigma_x$\\
\hline
\multirow{1}{*}{$\{y=0\}$}
& $\lambda\geq 1$ and $\mu\neq 0$ & y1& $\sigma_x$, $\sigma_0$\\ 
\hline
\multirow{1}{*}{$\{z=0\}$}
& $\mu=0$ and $b\geq 1$ & z1& $\sigma_x$ \\
\hline
\multirow{4}{*}{$\{xz=0\}$}
& $\mu\neq 0$ & xz1& $\sigma_0$ \\   
& $\mu=0$ and $b\geq 1$ & xz2 & $\sigma_x$\\
& $\mu=0, b<1$, and $a=\lambda$ & xz3 & $\sigma_x$\\
& $\mu \geq 1$ & xz4 & $\sigma_y$ \\
\hline
\multirow{3}{*}{$\{xy=0\}$}
&  & xy1 & $\sigma_0$ if $\lambda<1$ or $\mu=0$.\\  
&  & xy2 & $\sigma_x$ if $\lambda\geq 1$.\\  
&  & xy3 & $\sigma_y$ if $\mu\geq1$.\\  
\hline
\multirow{2}{*}{$\{yz=0\}$}
& & yz1& $\sigma_0$ if $\lambda<1$.\\
&  & yz2& $\sigma_0, \sigma_x$ if  $\lambda \geq 1$.\\
\hline
\multirow{3}{*}{$\{xyz=0\}$}
& & xyz1 & $\sigma_0$\\
& & xyz2 & $\sigma_x$ if $\lambda \geq 1$.\\ 
& & xyz3 & $\sigma_y$ if $\mu \geq 1$.\\
\hline
\end{tabular}
\caption{$\sigma_0=\{x=y=z=0\}$, $\sigma_x=\{x=z=0\}$ and $\sigma_y=\{ y=z=0\}$.}
\label{genposi}
\end{table}

\begin{example}\label{nasty}
Given  $\lambda>2$ and $0<b<1$, the surface $S$ with parametrization
$z=x^{\lambda}+x^{\lambda}y^b$
verifies the condition (xz3) of Theorem \ref{loglimitsteo}. 
Hence its logarithmic limits of tangents relatively to the divisor $\{xz=0\}$ is trivial. The proper inverse image of $S$ by the blow up with center $\{x=y=z=0\}$ admits the parametrization
$$\textstyle{\frac{z}{x}}=x^{\lambda-1}+x^{\lambda+b-1}\textstyle{\frac{y}{x}}^b.$$
By theorem \ref{loglimitsteo}, the logarithmic limit of tangents of $\widetilde{S}$ relatively to the divisor $\{x\textstyle{\frac{z}{x}}=0\}$ is not trivial.
\end{example}

\noindent
Example \ref{nasty} shows that the triviality of limits of tangents is not hereditary by blowing up. Lemma $\ref{nevermind}$ solves this problem.

\begin{lemma}\label{nevermind}
Let $N$ be the normal crossings divisor of a germ of  manifold $(M,o)$ of dimension three. Let $S$ be a quasi-ordinary surface of $M$ such that the logarithmic limit of tangents of $S$ along $N$ is trivial. Let $\pi:\widetilde{M}\rightarrow M$ be the blow up of $M$ along an admissible center for $S$ and $N$. Let $E$ be the exceptional divisor of $\pi$. Let $p \in \widetilde{S}\cap E$. If $S,N$ do not verify condition \em (xz3) \em of table \em (\ref{genposi}) \em at $o$,
$\widetilde{S}$ has trivial logarithmic limit of tangents along $\widetilde{N}$ at $p$
and $\widetilde{S}, \widetilde{N}$ do not verify condition \em (xz3) \em at $p$.
\end{lemma}

\begin{proof}
We will denote by (xy), (yz), (xyz) the situations (xyi), (yzi), (xyzi) for each i. We will assume that $b\neq +\infty$. The cases where $b= +\infty$ are much simpler. We also assume that after a blow up, the surface is not yet smooth.

\noindent
($\phi$1) We can assume that $z=x^{\lambda}y^{\mu}\varepsilon_1$.

\noindent
On the chart $(x, \textstyle{\frac{y}{x}}, \textstyle{\frac{z}{x}})$, $\widetilde{S}\cap \{x=0\}=\{\textstyle{\frac{z}{x}}=x=0\}$ and
$$\textstyle{\frac{z}{x}}=x^{\lambda+\mu-1}\textstyle{\frac{y}{x}}^{\mu}\varepsilon_2.$$
\noindent
Since $\lambda \geq \mu \geq 1$, we are in situation (x2) at each point of $\widetilde{S}\cap \{x=0\}$. The same happens in the chart $(\tfrac{x}{y},y,\tfrac{z}{y})$.

\noindent
($\phi$2) We can assume that $z=x^{\lambda}\delta_1+x^{a}y^{b}\varepsilon_1$. On the chart $(x,y,\tfrac{z}{x})$, $\widetilde{S}\cap \{x=0\}=\{\textstyle{\frac{z}{x}}-a=x=0\}$, $a=0$ if $\lambda>1$, $a\in \mathbb C^{\ast}$ otherwise, and

$$\tfrac{z}{x}=x^{\lambda-1}\delta_2+x^{a-1}y^{b}\varepsilon_2.$$

\noindent
If $\lambda>2$, we are in situation (x2) at each point of $\widetilde{S}\cap \{x=0\}$. Assume $1<\lambda<2$. Following a generalization of the proof of Theorem 3.5.5 of \cite{WALL}, $\widetilde{S}$ admits the parametrization
$$x=(\tfrac{z}{x})^{\frac{1}{\lambda-1}}\delta_3+(\tfrac{z}{x})^{\frac{a}{\lambda-1}-1}y^{b}\varepsilon_3.$$
Hence we are in situation (z1) at each point of $\widetilde{S}\cap \{x=0\}$. If $\lambda=1$, the situation is analogous to $\lambda>2$ or $1<\lambda<2$.

\noindent
(z1) The reasoning is similar to ($\phi$2) and we are in situation (xz2) or (x2) or (z1) at all points of $\widetilde{S}\cap E$.

\noindent
(xz1) Assume that $\lambda+\mu> 1$. On the chart $(x,\tfrac{y}{x},\tfrac{z}{x})$, $\widetilde{N}=\{x\tfrac{z}{x}=0\}$, $\widetilde{S}\cap \{x=0\}=\{x=\tfrac{z}{x}=0\}$ and
$$\tfrac{z}{x}=x^{\lambda+\mu-1}\tfrac{y}{x}^{\mu}\varepsilon_2.$$
\noindent
At the origin of the chart, if $\lambda>1$ we are in situation (xz1) or (xz4), otherwise we are in situation (yz) or (xz1) or (xz4). At a point of $\widetilde{S}\cap \{x=0\}$ where $\tfrac{y}{x}\neq 0$, the reasoning is similar to ($\phi$2) and we are in situation $(xz2)$. On the chart $(\tfrac{x}{y},y,\tfrac{z}{y})$, $\widetilde{N}=\{\tfrac{x}{y}y\tfrac{z}{y}=0\}$, $\widetilde{S}\cap \{y=0\}=\{y=\tfrac{z}{y}=0\}$ and
$$\tfrac{z}{y}=\tfrac{x}{y}^{\lambda}y^{\lambda+\mu-1}\varepsilon_2.$$
Hence, we are in situation (xyz) at the origin of the chart. At a point of $\widetilde{S}\cap \{y=0\}$ where $\tfrac{x}{y}\neq 0$, we are in situation $(xz2)$

\noindent
Assume that $\lambda+\mu=1$. On the chart $(x,\tfrac{y}{x},\tfrac{z}{x})$, after the change of coordinates $\widetilde{x}=\tfrac{z}{x}, \widetilde{y}=x, \widetilde{z}=\tfrac{y}{x}$, $\widetilde{N}=\{\widetilde{x}\widetilde{y}=0\}$ and
$$\widetilde{z}=\widetilde{x}^{\frac{1}{\mu}}\varepsilon_3.$$
We are in situation (xy2) at $(0,0,0)$. At a point of $\widetilde{S}\cap \{\widetilde{y}=0\}$ where
$\widetilde{x}\neq 0$, we are in situation (x2) or (z1). The reasoning for the chart $(\tfrac{x}{y},y,\tfrac{z}{y})$ is analogous.

\noindent
Assume that $\lambda+\mu<1$. On the chart $(x,\tfrac{y}{x},\tfrac{z}{x})$, $\widetilde{N}=\{x\tfrac{z}{x}=0\}$,
$$\tfrac{y}{x}=(\tfrac{z}{x})^{\frac{1}{\mu}}x^{\frac{1-\lambda-\mu}{\mu}}\varepsilon_2$$
and $\widetilde{S}\cap \{x=0\}=\{\tfrac{y}{x}=x=0\}$. At the origin of the chart we are in situation (xz1) or (xz4). At the remaining points of $\widetilde{S}\cap \{x=0\}$, the reasoning is similar to
the case ($\phi$2) and we are in situation (x2) or (z1). The reasoning on the chart $(\tfrac{x}{y},y,\tfrac{z}{y})$ is similar and the possible situations are (xyz) and (xz2).

\noindent
In the chart $(\tfrac{x}{z},\tfrac{y}{z},z)$, $\widetilde{N}=\{\tfrac{x}{z}z=0\}$ and
$$z=(\tfrac{x}{z})^{\frac{\lambda}{1-\lambda-\mu}}(\tfrac{y}{z})^{\frac{\mu}{1-\lambda-\mu}}\varepsilon_2.$$
\noindent
We are in the situation (xz1) or (xz4) or (yz) at $(0,0,0)$. Let $p$ be a point of $\widetilde{S}\cap \{z=0\}$. If $\tfrac{x}{z}\neq 0,\tfrac{y}{z}= 0$ at $p$, the reasoning is similar to ($\phi$2) and we are in situation (x2) or (z1). The reasoning is similar if $\tfrac{x}{z}= 0,\tfrac{y}{z}\neq 0$ at $o$ and we are in situation (xz2). If $\tfrac{x}{z},\tfrac{y}{z}\neq 0$ at $p$, $\widetilde{S}$ is smooth at $p$.

\noindent
(xz2) The reasoning is similar to ($\phi$2) and we are in situation (xz2) or (x2) or (z1) at all points of $\widetilde{S}\cap E$.

\noindent
(xz4) We can assume that $z=x^{\lambda}y^{\mu}\delta_1+x^ay^b\varepsilon_1$. On the chart $(x,y,\tfrac{z}{y})$, $\widetilde{N}=\{xy\tfrac{z}{y}=0\}$ and
$$\tfrac{z}{y}=x^{\lambda}y^{\mu-1}\delta_2+x^ay^{b-1}\varepsilon_2.$$
Hence we are in situation (xyz) at the origin. At the remaining points of $\widetilde{S}\cap \{y=0\}$
we are in situation (xz2), if $\mu>1$, and in situation (x2), if $\mu=1$.

\noindent
The remaining cases are similar to those studied in this proof.
\end{proof}

\begin{theorem}
Let $S$ be a quasi-ordinary surface of a germ of  manifold of dimension 3, $(M,o)$. Assume that  the limit of tangents of $S$ at $o$ is trivial. Let $M_0=M, \Gamma=\mathbb{P}_S^*M$. Let
$$M_0 \leftarrow M_1 \leftarrow M_2 \leftarrow \cdots \leftarrow M_m$$
be the sequence of blow ups that desingularizes $S$. Let $L_i$ be the center of the blow up $M_{i+1}\rightarrow M_i$ for $0 \leq i \leq m-1$. Let $S_i$ be the proper inverse image of $S$ by the map $M_i \rightarrow M_0$. Let $N_i$ be the inverse image of the first center by the map $M_i \rightarrow M_0$. Set $\Gamma_i=\mathbb{P}_{S_i}^*\langle M_i/N_i \rangle$, $\Lambda_i=\mathbb{P}_{L_i}^*\langle M_i/N_i \rangle$. Let $X_i$ be the blow up of  $\mathbb{P}^*\langle M_i/N_i \rangle$ along $\Lambda_i$. There are inclusion maps   $\mathbb{P}^*\langle M_{i+1}/N_{i+1} \rangle \hookrightarrow X_i$ such that the diagram (\ref{diafinal}) commutes.

\begin{equation}\label{diafinal}
\begin{array}{ccccccc}
  \mathbb{P}^*M_0 & \leftarrow &  \mathbb{P}^*\langle M_1/N_1 \rangle& \leftarrow & \cdots & \leftarrow &  \mathbb{P}^*\langle M_m/N_m \rangle  \\
\downarrow &  & \downarrow & & & & \downarrow\\  
M_0 & \leftarrow & M_1 & \leftarrow & \cdots & \leftarrow & M_m
\end{array}
\end{equation}
Moreover $\Gamma_m$ is a regular Lagrangean variety transversal to the set of poles of $\mathbb{P}\langle M_m/N_m \rangle$ and $\Gamma_m$ is the proper inverse image of $\Gamma_0$ by the map $\mathbb{P}^*\langle M_m / N_m \rangle$ $ \rightarrow \mathbb{P}^*M$.
\end{theorem}
\begin{proof}
It follows from \cite{mcewan},  Theorems \ref{inversecommutes}, \ref{diagteo} and \ref{loglimitsteo} and Lemma $\ref{nevermind}$.
\end{proof}

\bibliographystyle{amsplain}

\begin{thebibliography}{100} 



\bibitem{LIMITS} A. Ara\'ujo and O. Neto, \em Limits of Tangents of a Quasi-Ordinary Hypersurface, \em to be published at {Proc. Amer. Math. Soc.}



\bibitem{mcewan} 
{ C. Ban and L. Mcewan}, \em Canonical resolution of a quasi-ordinary surface singularity, \em  Canadian Journal of Mathematics, Vol. 52(6) (2000), 1149-1163.

\bibitem{milman} E. Bierstone and P. Milman, \em A simple constructive proof of canonical resolution of singularities, \em
Effective Methods in Algebraic Geometry, Progress in Math 94 (1991), 11-30.


\bibitem{kk} 
{ M. Kashiwara and T. Kawai,}  {\em
On holonomic systems of microdifferential equations. III: Systems
with regular singularities. \em} Publ. Res. Inst. Math. Sci. 17, (1981)
813-979.


\bibitem{ka3} { M. Kashiwara, Schapira, P.}, \em Sheaves on Manifolds\em, (Grundlehren Der Mathematischen Wissenschaften), Springer Verlag (2002).

\bibitem{LeTe1988} D.T. L\'e and B. Teissier, \em Limites d'espaces tangentes en g\'eometrie analitique\em, {Comment. Math. Helv.}, 63, (1988), 540-578.

\bibitem{Lipmanphd} J. Lipman, \em Quasi-ordinary singularities of embedded surfaces\em, {Ph.D. thesis, Harvard University}, 1965.


\bibitem{neto3} 
{ O. Neto},  \em Equisingularity and Legendrian curves\em, 
{Bull. London Math. Soc.} 33 (2001), 527-534.

\bibitem{neto4} 
{ O. Neto},  \em Blow up for a holonomic system. \em  Publ. Res. Inst. Math. Sci. 29 (1993), no. 2, 167-233.

\bibitem{neto5} 
{ O. Neto},  \em Systems of meromorphic microdifferential equations. \em Singularities and differential equations (Warsaw, 1993), 259-275, Banach Center Publ., 33, Polish Acad. Sci., Warsaw, 1996.


\bibitem{PHAM} F. Pham, \em Singularit\'{e}s des Syst\'{e}mes Differentiels de Gauss-Manin\em, (Progress in Mathematics), Birkhauser (1979).

\bibitem{WALL}, \em Singular points of Plane Curves\em, London math Society (2004)








\end{thebibliography}

\let\thefootnote\relax\footnotetext{\\ \\Jo\~ao Cabral, CMAF, Universidade de Lisboa Av. Gama Pinto, 2 1699-003 Lisboa, Portugal, and Departamento de Matem\'atica da Faculdade de Ci\^{e}ncias e Tecnologia da Universidade Nova de Lisboa, Quinta da Torre, 2829-516 Caparica, Portugal \\\email{joao.cabral.70@gmail.com}\\ \\
Orlando Neto, CMAF, Universidade de Lisboa Av. Gama Pinto, 2 1699-003 Lisboa, Portugal, and Faculdade de Ci\^{e}ncias da Universidade de Lisboa, Campo Grande 1749-016, Lisboa, Portugal \\\email{orlando60@gmail.com}\\ \\
Ant\'onio Ara\'ujo, DCeT, Universidade Aberta, R. Escola Polit\'ecnica 141-147, 1269-001 Lisboa, Portugal, and CMAF, Universidade de Lisboa Av. Gama Pinto, 2 1699-003 Lisboa, Portugal \\\email{ant.arj@gmail.com}}

\end{document}